

\def\W{{\cal W}}
\baselineskip=14pt
\parskip=10pt

\magnification=\magstephalf

\def\1{{\overline{1}}}
\def\2{{\overline{2}}}
\parindent=0pt
\overfullrule=0in

\def\frac#1#2{{#1 \over #2}}
\centerline
{\bf Increasing Consecutive Patterns in Words}
\bigskip
\centerline
{\it Mingjia YANG and Doron ZEILBERGER}

{\bf Abstract.} 
We show how to enumerate words in $1^{m_1} \dots n^{m_n}$ that avoid the increasing consecutive pattern $12 \dots r$ for any $r \geq 2$.  
Our approach yields an $O(n^{s+1})$ algorithm
to enumerate  words in $1^s \dots n^s$, avoiding the consecutive pattern $1\dots r$, for {\it any} $s$, and {\it any} $r$. 
This enables us to supply  many more terms to quite a few OEIS sequences, and create new ones.
We also treat the more general case of counting words with a specified number of the pattern of interest
(the avoiding case corresponding to zero appearances).
This article is accompanied by three Maple packages implementing our algorithms.

{\bf Introduction.} 

Rodica Simion and Herbert Wilf initiated the study of enumerating {\it classical} pattern-avoidance.
This is a very dynamic area with its own annual conference and Wikipedia page ([Wi]). 
Recall that 
a permutation $\pi=\pi_1 \dots \pi_n$ avoids a pattern $\sigma=\sigma_1 \dots \sigma_k$  if  none of the ${{n} \choose {k}}$ length-$k$ 
subsequences of  $\pi$, reduces to $\sigma$.

Alex Burstein ([Bu]), in a 1998 PhD thesis, under the direction of Herb Wilf, pioneered the enumeration of
{\it words} avoiding a set of patterns. This field is also fairly active today, with notable contributions
by, {\it inter alia}, Toufik Mansour ([BuM]) and Lara Pudwell ([P]).

The enumeration of permutations avoiding a given (classical) pattern, or a set of patterns, is notoriously difficult,
and it is widely believed to be intractable for most patterns, hence it would be nice to have other notions
for which the enumeration is more feasible. Such an analog was given, in 2003, by Sergi Elizalde
and Marc Noy, in a seminal paper ([EN]), that introduced the study of the enumeration of permutations
avoiding {\it consecutive} patterns. 
A permutation $\pi=\pi_1 \dots \pi_n$ avoids a consecutive pattern $\sigma=\sigma_1 \dots \sigma_k$  
if  none of the $n-k+1$ length-$k$ consecutive subwords, $\pi_i \pi_{i+1} \dots \pi_{i+k-1}$ of  $\pi$, reduces to $\sigma$.

Algorithmic approaches to the enumeration of {\it permutations} avoiding sets of consecutive patterns were given
by Brian Nakamura, Andrew Baxter, and Doron Zeilberger ([Na], [BaNaZ]). Our present approach 
may be viewed as an extension, from permutations to words, of Nakamura's paper, who was also inspired by
the Goulden-Jackson cluster method, but in a sense, is more straightforward, and closer in spirit to the original
Goulden-Jackson cluster method ([GJ], that is beautifully exposited (and extended!) in [NoZ]). 

In this article we will consider consective patterns of the form $1 \dots r$, i.e. {\it increasing consecutive patterns},
and show how to count words in $1^{m_1} \dots n^{m_n}$ avoiding the pattern  $1 \dots r$ (Theorem 1, that is due to Ira Gessel[Ge]). 
Throughout this article we will only consider consecutive patterns, so the word ``consecutive'' may be omitted. 
In particular, we will look at how to efficiently count words in $1^{s} \dots n^{s}$ avoiding the pattern $1 \dots r$.
All the sequences for $s=1$ and $3 \leq r \leq 9$ are  in the On-Line Encyclopedia of Integer Sequences, with many
terms.
Also, quite a few of theses sequences for $s>1$ are already there, but with very few terms.
Our implied algorithms are $O(n^{s+1})$ and hence yield many more terms, and, of course, new sequences.

In the last part of the paper, we will provide a new proof of Theorem 1 by tweaking the Goulden-Jackson cluster method. Using this proof, along with a little more effort, we will generalize Theorem 1 to counting words with a specified number of the pattern 
$12 \dots r$ (Theorem 2), instead of just {\it avoiding}, that is, having {\it zero} occurrence of the pattern of interest. 

We close this introduction by mentioning the pioneering work of Anthony Mendes and Jeff Remmel([MR]), 
in combining the two keywords `consecutive patterns' and `words'. We were greatly inspired by their article,
but our focus is algorithmic.

{\bf Maple Packages}: This article is accompanied by three Maple packages available from the webpage:

{\tt http://sites.math.rutgers.edu/\~{}zeilberg/mamarim/mamarimhtml/icpw.html}  \quad .

These are

$\bullet$ {\tt ICPW.txt}: For fast  enumeration of sequences enumerating words avoiding increasing consecutive patterns.

$\bullet$ {\tt ICPWt.txt}: For fast  computation of  sequences of weight-enumerators for  words according to the number of 
increasing consecutive patterns ($t=0$ reduces to the former case).

$\bullet$ {\tt GJpats.txt}: For conjecturing generating functions (that still have to be proved by humans).

This page also has links to numerous input and output files. The input files can be modified to generate more data, if desired.

{\bf The Goulden-Jackson Cluster Method}

Recall that the original Goulden-Jackson method ([GJ][NoZ]) inputs a {\it finite} alphabet, $A$, 
that may be taken to be $\{1, ..., n\}$, and a finite set of `bad words', $B$.

It outputs a certain {\bf rational function}, let's call it $F(x_1, \dots, x_n)$, that is the multi-variable generating function,
in $x_1, \dots, x_n$,  for the discrete $n$-variable function
$$
f(m_1, \dots, m_n) \quad ,
$$
that counts the words in $1^{m_1} \dots n^{m_n}$ (there are altogether $(m_1 + \dots + m_n)!/(m_1! \cdots m_n!)$ of them)
that {\bf never} contain, as {\it consecutive} subwords (aka {\it factors} in linguistics) any member of $B$.
In other words:
$$
F(x_1, \dots , x_n)= \sum_{(m_1, \dots, m_n) \in N^n}  f(m_1, \dots, m_n) \, x_1^{m_1} \cdots x_n^{m_n} \quad .
$$
This is nicely implemented in the Maple package {\tt DavidIan.txt}, that accompanies [NoZ], and is freely available from

{\tt http://sites.math.rutgers.edu/\~{}zeilberg/tokhniot/DavidIan.txt } \quad .

For example, if $n=4$, so the alphabet is $\{1,2,3,4\}$ and the set of `bad words' to avoid is $\{1234,1432\}$, then,
starting a Maple session, and typing:

{\tt read `DavidIan.txt`: lprint(subs(t=0,GJgf({1,2,3,4},{[1,2,3,4],[1,4,3,2]},x,t)));}

immediately returns

{\tt 1/(1-x[1]-x[2]-x[3]-x[4]+ 2*x[1]*x[2]*x[3]*x[4])} \quad ,

that in Humanese reads
$$
\frac{1}{1-x_1-x_2-x_3-x_4 \, + \, 2\, x_1 x_2 x_3 x_4} \quad .
$$

{\bf Enumerating Words Avoiding Consecutive Patterns: Let the Computer Do the Guessing}

{\it Now} we are interested in words in an {\it arbitrarily large} alphabet $\{1, \dots , n\}$ avoiding a set of consecutive patterns,
but each pattern, e.g. $123$, entails  an {\it arbitrarily large} set of forbidden consecutive subwords. For example, in this case, the set of forbidden consecutive subwords is 
$$
\{ i_1 \, i_2 \, i_3 \, |\, 1 \leq i_1 < i_2 < i_3 \leq n \} \quad .
$$
We can ask {\tt DavidIan.txt} to find the generating function for each specific $n$, and then hope to conjecture
a {\bf general} formula in terms of $x_1, \dots , x_n$, for {\it general} (i.e. {\bf symbolic}) $n$.

This is accomplished by the Maple package {\tt GJpats.txt}, available from the webpage of this article.
It uses the original  {\tt DavidIan.txt} to produce the corresponding generating functions for increasing values for $n$, and then attempts to conjecture a
{\it meta-pattern}.
For example for words avoiding the consecutive pattern $123$ (alias the word $123$), for $n=3$,

{\tt GFpats($\{[1,2,3]\},x,3,0$);} yields
$$
1/(1-x_1-x_2-x_3+ x_1x_2x_3) \quad .
$$
This is simple enough. Moving right along,

{\tt GFpats($\{[1,2,3]\},x,4,0$);} yields
$$
1/(1-x_1-x_2-x_3-x_4+ x_1x_2x_3+ x_1x_2x_4+ x_1x_3x_4+ x_2x_3x_4 -x_1x_2x_3x_4) \quad ,
$$
while {\tt GFpats($\{[1,2,3]\},x,5,0$);} yields 
$$
1/(1-x_1-x_2-x_3-x_4 -x_5+ x_1x_2x_3+ x_1x_2x_4+ x_1x_2x_5+x_1x_3x_4+ x_1 x_3 x_5 + x_1x_4x_5+
$$
$$
 x_2x_3x_4 +x_2x_3x_5+x_2x_4x_5+x_3x_4x_5 -x_1x_2x_3x_4 -x_1x_2x_3x_5 - x_1x_2x_4x_5 - x_1x_3x_4x_5 -x_2x_3x_4x_5) \quad .
$$

These look like {\it symmetric} functions. Procedure {\tt  SPtoM(P,x,n,m)} expresses a polynomial, P, in the indexed variables
$x[1], \dots, x[n]$ in terms of the {\it monomial symmetric polynomials} $m_\lambda$. Applying this procedure we have

{\tt SPtoM(denom(GFpats($\{[1,2,3]\},x,5,0)),x,5,m$);} yields

{\tt  -m[1, 1, 1, 1] + m[1, 1, 1] - m[1] + m[]} \quad .

{\tt SPtoM(denom(GFpats($\{[1,2,3]\},x,6,0)),x,6,m$);} yields

{\tt m[1,1,1,1,1,1]-m[1,1,1,1]+m[1,1,1]-m[1]+m[]} \quad .

{\tt SPtoM(denom(GFpats($\{[1,2,3]\},x,7,0)),x,7,m$);} yields

{\tt -m[1,1,1,1,1,1,1]+m[1,1,1,1,1,1]-m[1,1,1,1]+m[1,1,1]-m[1]+m[]} \quad .

You don't have to be a Ramanujan to conjecture the following result.

{\bf Fact:} The generating function for words in $\{1, 2, \dots , n\}$ avoiding the consecutive pattern $123$, let's call it 
$F_3(x_1,\dots, x_n)$  is
$$
F_3(x_1, \dots , x_n) = \frac{1}{1-e_1+e_3-e_4+e_6-e_7+ e_9-e_{10} \, + \dots } \quad ,
$$
where $e_i$ stands for the {\it elementary symmetric function}  of degree $i$ in $x_1, \dots, x_n$, i.e. the coefficient
of $z^i$ in $(1+x_1\,z)\dots (1+ x_n\, z)$.

(Note that $e_i=m_{1^i}$).

Doing the analogous guessing for the consecutive patterns $1234$ and $12345$, a {\it meta-pattern} emerges, and we
were safe in formulating the following theorem that we discovered using the present experimental mathematics approach.
After the first version of this article was posted, we found out, thanks to Justin Troyka, that
this theorem is due to Ira Gessel ([Ge] , p. 51, Example 3).

{\bf Theorem 1 (Gessel [Ge])} For $n \geq 1, r\geq 2$,
the generating function for words in $\{1, 2, \dots , n\}$ avoiding the consecutive pattern $12\dots r$, let's call it 
$F_r(x_1,\dots, x_n)$  is
$$
F_r(x_1, \dots , x_n) = \frac{1}{1-e_1+e_r-e_{r+1}+e_{2r}-e_{2r+1}+e_{3r}-e_{3r+1}+  \dots } \quad .
$$

Of course, if Gessel did not prove it before us,  these would have been `only' guesses, but once known, humans can prove them.
We did it  by tweaking the cluster method to apply to an {\it arbitrarily large} alphabet, i.e. where even the size of the alphabet, $n$, is {\it symbolic}.
Our proof of Gessel's theorem will be given at the end of this article.

{\bf Efficient Computations}

The Theorem immediately implies the following partial recurrence equation for the actual coefficients.

\vfill\eject

{\bf Fundamental Recurrence}: Let $f_r({\bf m})$ be the number of words in the alphabet $\{1, \dots , n \}$ with
$m_1$ $1$'s, $m_2$ $2$'s, \dots ,  $m_n$ $n$'s (where we abbreviate ${\bf m}=(m_1, \dots, m_n)$) that avoid the consecutive pattern $1\dots r$. Also
let $V_i$ be the set of $0-1$ vectors of length $n$ with $i$ ones, then

$$
f_r({\bf m})= \sum_{ {\bf v} \in V_1} f_r({\bf m} - {\bf v}) \, - \,\sum_{ {\bf v} \in V_r} f_r({\bf m} - {\bf v}) 
$$
$$
+\sum_{ {\bf v} \in V_{r+1}} f_r({\bf m} - {\bf v}) \, - \,\sum_{ {\bf v} \in V_{2r}} f_r({\bf m} - {\bf v}) 
$$
$$
+\sum_{ {\bf v} \in V_{2r+1}} f_r({\bf m} - {\bf v}) \, - \,\sum_{ {\bf v} \in V_{3r}} f_r({\bf m} - {\bf v})  
$$
$$
+\sum_{ {\bf v} \in V_{3r+1}} f_r({\bf m} - {\bf v}) \, - \,\sum_{ {\bf v} \in V_{4r}} f_r({\bf m} - {\bf v})  +\, \dots \quad  .
$$

Suppose that we want to compute $f_3(1^{100})$, i.e. the number of permutations of length $100$ that avoid the consecutive
pattern $123$. If we use the above recurrence literally, we would need about $2^{100}$ computations, but there is
a shortcut!

{\bf Enter Symmetry}

It follows from the symmetry of the generating function $F_r(x_1, \dots, x_n)$, that $f_r(m_1, \dots, m_n)$ is symmetric, hence the above Fundamental Recurrence 
immediately implies the following recurrence, that enables  a very fast computation of the sequences, let's call them
$a_r(n)$, for the number of {\it permutations} of length $n$ that avoid the consecutive pattern $1\dots r$.

{\bf Fast Recurrence For Enumerating Permutations avoiding the consecutive pattern $1 \dots r$}

$$
a_r(n)= n a_r(n-1) -{{n} \choose {r}} a_r(n-r) + {{n} \choose {r+1}} a_r(n-r-1) -{{n} \choose {2r}} a_r(n-2r) + {{n} \choose {2r+1}} a_r(n-2r-1)  
$$
$$
-{{n} \choose {3r}} a_r(n-3r)   +  {{n} \choose {3r+1}} a_r(n-3r-1)   \, - \dots \quad .
$$
This recurrence  goes back to Florence Nightingale David and David Barton ([DB], p. 157, line 6 from the top),
whose proof used a probabilistic language and an inclusion-exclusion argument that may be viewed as a precursor of the cluster method, 
applied to the special case of increasing patterns.

Equivalently, we have the following exponential generating function ([DB], p. 157, line 4)
$$
\sum_{n=0}^{\infty} a_r(n) \frac{x^n}{n!} \, = \,
\frac{1}{1-x+ \frac{x^r}{r!} - \frac{x^{r+1}}{(r+1)!} + \frac{x^{2r}}{(2r)!}   
-\frac{x^{2r+1}}{(2r+1)!}+ \frac{x^{3r}}{(3r)!} -  \frac{x^{3r+1}}{(3r+1)!} + \dots } \quad .
$$
While this `explicit' (exponential) generating function is `nice', it is more  efficient to use the fast recurrence. And
indeed, the OEIS has these sequences for $3 \leq r \leq 9$, with many terms.
These are (in order): {\tt A049774, A117158, A177523, A177533, A177553, A230051, A230231}.

{\bf Efficient Computations of Permutations of words with Two Occurrences of each Letter}

Let $b_r(n)$ be the number of words with $2$ occurrences of each of $1, 2, \dots, n$ avoiding the pattern $1\dots r$.
Quite a few of them are currently (April 17, 2018) in the OEIS, but with relatively few terms

$\bullet$ $b_3(n)$: {\tt https://oeis.org/A177555} (15 terms)

$\bullet$ $b_4(n)$: 
{\tt https://oeis.org/A177558}  (15 terms)

$\bullet$ $b_5(n)$: 
{\tt https://oeis.org/A177564} (14 terms)

$\bullet$ $b_6(n)$: 
{\tt https://oeis.org/A177574} (14 terms)

$\bullet$ $b_7(n)$:
{\tt https://oeis.org/A177594}  (14 terms)

$b_r(n)$ for $r>7$ are not yet (April 17, 2018) in the OEIS.

We can compute $b_r(n)$ in cubic time as follows. If you plug-in $f_r(2^n)$ into the Fundamental Recurrence, you are
forced to consider the more general quantities of the form $f_r(2^\alpha 1^\beta)$. Defining
$$
B_r(\alpha,\beta)=f_r(2^\alpha 1^\beta) \quad ,
$$
and using symmetry, we get the following recurrence for $B_r(\alpha,\beta)$.
$$
B_r(\alpha,\beta)=
\alpha B_r(\alpha-1,\beta+1) + \beta B_r(\alpha,\beta-1)
$$
$$
-\sum_{i_1+i_2=r} {{\alpha} \choose {i_1}} {{\beta} \choose {i_2}} B_r(\alpha-i_1,\beta-i_2+i_1)
+\sum_{i_1+i_2=r+1} {{\alpha} \choose {i_1}} {{\beta} \choose {i_2}} B_r(\alpha-i_1,\beta-i_2+i_1)
$$
$$
-\sum_{i_1+i_2=2r} {{\alpha} \choose {i_1}} {{\beta} \choose {i_2}} B_r(\alpha-i_1,\beta-i_2+i_1)
+\sum_{i_1+i_2=2r+1} {{\alpha} \choose {i_1}} {{\beta} \choose {i_2}} B_r(\alpha-i_1,\beta-i_2+i_1)\,  - \,\dots  \quad .
$$
In particular $b_r(n)=B_r(n,0)$. Using this recurrence we (easily!) obtained 80 terms of each of  the sequences
$b_r(n)$ for $3 \leq r \leq 9$, and could get many more. See the output file

{\tt http://sites.math.rutgers.edu/\~{}zeilberg/tokhniot/oICPW1.txt} \quad .

{\bf Efficient Computations of Permutations of words with Three Occurrences of each Letter}

Let $c_r(n)$ be the number of words with $3$ occurrences of each of $1, 2, \dots, n$ avoiding the pattern $1\dots r$.
Quite a few of them are currently (April 17, 2018) in the OEIS, but with relatively few terms

$\bullet$ $c_3(n)$: {\tt https://oeis.org/A177596} (Only 10 terms)

$\bullet$ $c_4(n)$: {\tt https://oeis.org/A177599} (Only 10 terms)

$\bullet$ $c_5(n)$: {\tt https://oeis.org/A177605} (Only 10 terms)

$\bullet$ $c_6(n)$: {\tt https://oeis.org/A177615} (Only 9 terms)

$\bullet$ $c_7(n)$: {\tt https://oeis.org/A177635} (Only 9 terms)

$c_r(n)$ for $r>7$ are not yet in the OEIS.

We can compute $c_r(n)$ in quartic time as follows. If you plug-in $f_r(3^n)$ into the Fundamental Recurrence, you are
forced to consider the more general quantities of the form $f_r(3^\alpha 2^\beta 1^\gamma)$. Defining
$$
C_r(\alpha,\beta, \gamma)=f_r(3^\alpha 2^\beta 1^\gamma ) \quad ,
$$
and using symmetry, we get the following recurrence for $C_r(\alpha,\beta, \gamma)$.
$$
C_r(\alpha,\beta, \gamma)=
\alpha C_r(\alpha-1,\beta+1, \gamma) + \beta C_r(\alpha,\beta-1, \gamma+1) + \gamma C_r(\alpha,\beta, \gamma-1)
$$
$$
- \, \sum_{i_1+i_2+i_3=r} {{\alpha} \choose {i_1}} {{\beta} \choose {i_2}}  {{\gamma} \choose {i_3}} C_r(\alpha-i_1,\beta-i_2+i_1, \gamma-i_3+i_2)
$$
$$
+ \, \sum_{i_1+i_2+i_3=r+1} {{\alpha} \choose {i_1}} {{\beta} \choose {i_2}}  {{\gamma} \choose {i_3}} C_r(\alpha-i_1,\beta-i_2+i_1, \gamma-i_3+i_2)
$$
$$
-\sum_{i_1+i_2+i_3=2r} {{\alpha} \choose {i_1}} {{\beta} \choose {i_2}}  {{\gamma} \choose {i_3}} C_r(\alpha-i_1,\beta-i_2+i_1, \gamma-i_3+i_2)
$$
$$
+ \, \sum_{i_1+i_2+i_3=2r+1} {{\alpha} \choose {i_1}} {{\beta} \choose {i_2}}  {{\gamma} \choose {i_3}} C_r(\alpha-i_1,\beta-i_2+i_1, \gamma-i_3+i_2) \, - \, \dots
$$
In particular, $c_r(n)=C_r(n,0,0)$. Using this recurrence we (easily!) obtained 40 terms of each of  the sequences
$c_r(n)$ for $3 \leq r \leq 9$, and could get many more. See the output file

{\tt http://sites.math.rutgers.edu/\~{}zeilberg/tokhniot/oICPW1.txt} \quad .

{\bf Efficient Computations of Permutations of words with Four Occurrences of each Letter}

Let $d_r(n)$ be the number of words with $4$ occurrences of each of $1, 2, \dots, n$ avoiding the pattern $1\dots r$.
Quite a few of them are currently (April 17, 2018) in the OEIS, but with relatively few terms.

$\bullet$ $d_3(n)$: {\tt https://oeis.org/A177637} (8 terms)

$\bullet$ $d_4(n)$:  {\tt https://oeis.org/A177640} (8 terms)

$\bullet$ $d_5(n)$: {\tt https://oeis.org/A177646} (8 terms)

$\bullet$ $d_6(n)$: {\tt https://oeis.org/A177656} (8 terms)

$\bullet$ $d_7(n)$: {\tt https://oeis.org/A177676} (8 terms)

$d_r(n)$ for $r>7$ are not yet in the OEIS.

We can compute $d_r(n)$ in quintic time as follows. If you plug-in $f_r(4^n)$ into the Fundamental Recurrence, you are
forced to consider the more general quantities of the form $f_r(4^\alpha 3^\beta 2^\gamma 1^\delta)$. Defining
$$
D_r(\alpha,\beta, \gamma, \delta)=f_r(4^\alpha 3^\beta 2^\gamma 1^\delta) \quad ,
$$
and using symmetry, we get the following recurrence for $D_r(\alpha,\beta, \gamma, \delta)$.
$$
D_r(\alpha,\beta, \gamma, \delta)=
\alpha D_r(\alpha-1,\beta+1, \gamma, \delta) + \beta D_r(\alpha,\beta-1, \gamma+1, \delta) + \gamma D_r(\alpha,\beta, \gamma-1,\delta+1)
+ \delta D_r(\alpha,\beta, \gamma,\delta-1)
$$
$$
- \sum_{i_1+i_2+i_3+i_4=r} {{\alpha} \choose {i_1}} {{\beta} \choose {i_2}}  {{\gamma} \choose {i_3}}  {{\delta} \choose {i_4}} 
D_r(\alpha-i_1,\beta-i_2+i_1, \gamma-i_3+i_2, \delta-i_4+i_3)
$$
$$
+\sum_{i_1+i_2+i_3+i_4=r+1} {{\alpha} \choose {i_1}} {{\beta} \choose {i_2}}  {{\gamma} \choose {i_3}}  {{\delta} \choose {i_4}} 
D_r(\alpha-i_1,\beta-i_2+i_1, \gamma-i_3+i_2, \delta-i_4+i_3)
$$
$$
-\sum_{i_1+i_2+i_3+i_4=2r} {{\alpha} \choose {i_1}} {{\beta} \choose {i_2}}  {{\gamma} \choose {i_3}}  {{\delta} \choose {i_4}} 
D_r(\alpha-i_1,\beta-i_2+i_1, \gamma-i_3+i_2, \delta -i_4+i_3) 
$$
$$
+\sum_{i_1+i_2+i_3+i_4=2r+1} {{\alpha} \choose {i_1}} {{\beta} \choose {i_2}}  {{\gamma} \choose {i_3}}  {{\delta} \choose {i_4}} 
D_r(\alpha-i_1,\beta-i_2+i_1, \gamma-i_3+i_2, \delta -i_4+i_3) \, -\, \dots 
$$

In particular $d_r(n)=D_r(n,0,0,0)$. Using this recurrence we (easily!) obtained 20 terms of each of  the sequences
$c_d(n)$ for $3 \leq r \leq 9$, and could get many more. See the output file

{\tt http://sites.math.rutgers.edu/\~{}zeilberg/tokhniot/oICPW1.txt} \quad .

{\bf Keeping Track of the Number of Occurrences}

Above we showed how to enumerate words {\it avoiding} the consecutive pattern $1 \dots r$, in other words,
the number of words, with a specified number of each letters, with {\bf zero} such patterns.
 With a little more effort we can answer the more general question about the number of such words with
exactly $j$ consecutive patterns $1 \dots r$ for {\it any} $j$, not just $j=0$. 
Let $\W({\bf m})=\W(m_1, \dots, m_n)$ be the set of words in the alphabet $1, \dots , n$ with $m_1$ $1$'s, $\dots$, $m_n$ $n$'s
(note that the number of elements of $\W({\bf m})$ is $(m_1 + \dots +m_n)!/(m_1! \cdots m_n!)$).

We are interested in the polynomials in $t$
$$
g_r({\bf m}; t)=\sum_{w \in \W({\bf m}) }  t^{\alpha(w)} \quad ,
$$
where $\alpha(w)$ is the number of occurrences of the consecutive pattern  $1\dots r$ in the word $w$.
(For example $\alpha(831456178)=3$ if $r=3$. Note that $\alpha(w)=0$ if $w$ avoids the pattern.)

[Also note that $g_r({\bf m}; 0)=f_r({\bf m})$  and $g_r({\bf m}; 1)=(m_1 + \dots + m_n)!/(m_1! \cdots m_n!)$.]

Using GJpats.txt we were able to conjecture the following theorem, whose proof will be presented later.

We first need to define certain families of  polynomial sequences.

{\bf Definition}: For any integer $k \geq 1$ and $r\geq 2$,  $P_k^{(r)}(t)$ is defined as follows.

If $k<r$, then it is $0$. If $k=r$ then it is $t-1$, and if $k>r$ then
$$
P_k^{(r)}(t) \, = \, (t-1) \sum_{i=1}^{r-1} P_{k-i}^{(r)}(t) \quad .
$$

{\bf Theorem 2:} For $k \geq 1, r\geq 2$,
the generating function of $g_r({\bf m};t)$, let's call it $G_r(x_1, \dots, x_n;t)$, is
$$
G_r(x_1, \dots , x_n;t) = \frac{1}{1-e_1- \sum_{k=r}^{n} P_k^{(r)}(t) e_k} \quad .
$$

This implies the

{\bf Fundamental Recurrence for $g_r$}: Let $g_r({\bf m};t)$ be the weight-enumerator of words in the alphabet $\{1, \dots , n \}$ with
$m_1$ $1$'s, $m_2$ $2$'s, \dots $m_n$ $n$'s (where we abbreviate ${\bf m}=(m_1, \dots, m_n)$), according to the weight 

``$t$ raised to the power of the number of occurrences of the consecutive pattern $1 \dots r$''. 

Also, let $V_k$ be the set of $0-1$ vectors of length $n$ with $k$ ones, then

$$
g_r({\bf m})= \sum_{ {\bf v} \in V_1} g_r({\bf m} -{\bf v}) \, + \, \sum_{k=r}^{n} \sum_{ {\bf v} \in V_k} P_k^{(r)}(t)\, g_r({\bf m} -{\bf v}) \quad .
$$

Analogously to the avoidance case we can get efficient recurrences for the permutations, and words in $1^s \cdots n^s$,  for each $s>1$. For each $s$ it is still
polynomial time, but things are slower because of the variable $t$. This is implemented in the Maple package {\tt ICPWt.txt} \quad .

{\bf Proofs.}

{\bf A New Proof of Gessel's Theorem 1.}

We will use the general set-up of the Goulden-Jackson cluster method as described in [NoZ], but 
will be able to make things simpler by taking advantage of the specific structure of our forbidden patterns,
that happen to be the increasing patterns $1 \dots r$.
That would enable us to use an elegant combinatorial argument, without solving a system of linear equations.

First let us quickly review some basic definitions. (We will not go into the details of the cluster method but readers 
who wish to see an excellent and concise summary of the cluster method are welcome to refer to the first section of [W].) 
A {\bf marked word} is a word with some of its factors (consecutive subwords) marked. We are assuming that all the marks are 
in the set of bad words $B$. For example (13212; [1,3]) is a marked word with 132 marked, with [1,3] denoting the location of the mark. 
A {\bf cluster} is a marked word where the adjacent marks overlap with each other and all the letters in the 
underlying word belong to at least one mark of the cluster.
For example (145632; [1,3],[2,4],[4,6]) is a cluster whereas (145632; [1,3],[4,6]) is not. We let the weight of a marked word 
$w=w_1w_2 \dots w_k$ be $weight(w):=(-1)^{|S|} \cdot \prod_{i=1}^{k}{x[w_i]}$ where $S$ is the set of marks in $w$. 
For example, the weight of the cluster (135632; [1,3],[2,4],[4,6]) is $(-1)^3 x_1x_2x_3^2x_5x_6$.

Let $M$ be the set of all marked words in the alphabet $\{1,..,n\}$. Recall from [NoZ] that 
$weight(M) \, = \, 1+ weight(M)\cdot (x_1+x_2+\dots+x_n)$+$weight(M)\cdot weight(C)$ 
where $C$ is the set of all possible clusters. This implies, according to [NoZ], 
that the multivariate generating function for words avoiding the consecutive pattern $1 \dots r$ (i.e. our target generating function) 
is equal to $weight(M)=\frac{1}{1-e_1-weight(C)}$. So we only need to figure out 
$weight(C)$. However, to use the classical Goulden-Jackson cluster method, we would have 
to solve a system of ${n}\choose{r}$ (the number of bad words) equations (recall that we write C as a summation of 
$C[v]$'s where $v$ is a word in B, and for each C[v] there is an equation) and no obvious symmetry argument seems to help.
So we will use a slick combinatorial approach.

Notice that since the pattern to be avoided is $12 \dots r$, the clusters can only be of the form 
$$(a_1 \dots a_j;[1,r],\dots)$$ where $$1 \leq a_1<a_2<\dots<a_j \leq n \quad. $$ 
Therefore $weight(C)$ is a summations of multivariate monomials on $x_1,x_2,..,x_n$ where the exponent of each variable $x_i$ is zero or one.

Any fixed monomial in $weight(C)$, it can come from many different clusters. The number of clusters it comes from and the coefficient 
of the monomial are uniquely determined by the number of variables in the monomial. For example, for $r=3$, the monomial $x_1x_3x_5x_6x_7$ 
can come from the cluster $(13567;[1,3],[2,4],[3,5])$ or $(13567;[1,3],[3,5])$. The first cluster contributes weight $(-1)^3x_1x_3x_5x_6x_7$ 
whereas the second cluster contributes weight  $(-1)^2x_1x_3x_5x_6x_7$. So when summing up, they cancel each other out and there is no monomial $x_1x_3x_5x_6x_7$ in $weight(C)$. 
So is the case with any other monomial of degree $5$.
Therefore, let us focus on the monomial $x_1x_2x_3 \dots x_k$ and figure out its coefficient. 

{\bf{Definition}}: Let coeff$(x_1x_2 \dots x_k)$ ($k \geq 1 $) be the  coefficient of $x_1x_2 \dots x_k$ in $weight(C)$.

It is clear that for $k<r$, coeff($x_1x_2x_3 \dots x_k)=0$, because $12\dots k$ cannot be a cluster (it does not have enough letters to be marked). 
And when $k=r$, we have coeff$(x_1x_2 \dots x_k)=-1$, since there can be only one mark. So let us move on to the case when $k>r$. 
We have the following claim.

{\bf{Claim 1:}} 

For $k>r$, coeff$(x_1x_2 \dots x_k)=-$ coeff$(x_2x_3 \dots x_k)-$ coeff$(x_3x_4 \dots x_k)-\dots-$ coeff$(x_rx_{r+1} \dots x_k)$. 
(i.e. coeff$(x_1x_2 \dots x_k)=-$ coeff$(x_1x_2 \dots x_{k-1})-$ coeff$(x_1x_2 \dots x_{k-2})- \dots -$ coeff$(x_{1}x_{2} \dots x_{k-r+1})$.)

This is because there are $(r-1)$ ways
in which the left-most marked word can `interface' with the one to its immediate right.
For example, if the clusters are of the form $(1\dots k;[1,r],[3,r+2],\dots)$ (that is, the 
second mark starts at 3), then the contribution  will be $(-1) \cdot$ coeff$(x_3x_4 \dots x_k)$. 
This is simply because of the bijection between the set of clusters in the form of $(1\dots k;[1,r],[3,r+2],\dots)$ with set of the 
clusters in the form $(3 \dots k;[3,r+2],
\dots)$. By peeling off the first mark $[1,r]$, we just lose a factor of $(-1)$ in the coefficient of our monomial. 

Similarly, if the clusters are of the form $(1 \dots k;[1,r],[u,u+r-1],\dots)$ ($1<u \leq r$), then the contribution from this case
will be $(-1) \cdot$ coeff$(x_ux_{u+1} \dots x_k)$. Note that if $k<2r-1$, there cannot be as many as $(r-1)$ cases.
However, in this case, we can make the convention that there are $(r-1)$ places for the second mark because 
for $k<r$ the coefficient of $x_1x_2x_3 \dots x_k$ is 0. So the above formula still holds. For example, for the clusters associated with the 
word 123456, and $r=4$, the first mark has to be 1234, the second mark can only be 2345 or 3456. But, according to the
natural convention, the second mark can also 
start with 4 and be 456, 
and so, coeff$(x_1x_2x_3x_4x_5x_6)=-$coeff$(x_2x_3x_4x_5x_6)-$coeff$(x_3x_4x_5x_6)-$coeff$(x_4x_5x_6)$$=-$coeff$(x_2x_3x_4x_5x_6)-$coeff$(x_3x_4x_5x_6)$.

So we have:
coeff$(x_1x_2 \dots x_r)=-1$; coeff$(x_1x_2 \dots x_{r+1})=(-1) \cdot (-1)=1$; 
coeff$(x_1x_2 \dots x_{r+2})= -$coeff$(x_2x_3 \dots x_{r+2})$ $-$ coeff$(x_3x_4 \dots x_{r+2})$
$= -$coeff$(x_1x_2 \dots x_{r+1})$ $-$ coeff$(x_1x_2 \dots x_{r})=0$. Continuing this process, it is easy to see that 
$x_1x_2 \dots x_{mr}$ $(m \geq 1)$ has coefficient $-1$ (so is any other monomial of degree $mr$ ) and $x_1x_2 \dots x_{mr+1}$ has coefficient $1$ 
(so is any other monomial of degree $mr+1$). The monomials with other number of variables all have coefficient 0. From this argument and 
summing over all clusters, we conclude $weight(C)=-e_r+e_{r+1}-e_{2r}+e_{2r+1}+ \dots$ and therefore 
$weight(M)=\frac{1}{1-e_1+e_r-e_{r+1}+e_{2r}-e_{2r+1}+ \dots}$.

{\bf Proof of Theorem 2.}

This proof can be directly generalized from the proof of Theorem 1 based on the `$t$-generalization' described in [NoZ]. 
Again, let the set of marked words on $\{1,2,\dots,n\}$ be $M$. However, this time we let the weight of a marked word $w$ of length $k$ 
be $weight(w):=(t-1)^{|S|} \cdot \prod_{i=1}^{k}{x[w_i]}$ where $S$ is the set of marks in $w$. We still have 
$weight(M)\, = \, 1+ weight(M)\cdot (x_1+x_2+\dots+x_n)$+$weight(M)\cdot weight(C)$ and $G_r(x_1, \dots , x_n;t)$ is equal to $weight(M)$, 
which is $\frac{1}{1-e_1-weight(C)}$.

The procedure to calculate $weight(C)$ directly follows from the proof of Theorem 1. We simply replace $(-1)$ with $(t-1)$ in various places, 
because the only difference is that now we assign a different weight to a marked word. For example, we have coeff$(x_1x_2 \dots x_r)=t-1$; 
coeff$(x_1x_2 \dots x_{r+1})=(t-1)(t-1)=(t-1)^2$; coeff$(x_1x_2 \dots x_{r+2})= (t-1)$(coeff$(x_2x_3 \dots x_{r+2})$ $+
$ coeff$(x_3x_4 \dots x_{r+2}))$ $=(t-1)((t-1)+(t-1)^2).$ Again it is clear that for $k<r$, coeff($x_1x_2x_3 \dots x_k)=0$ and when $k=r$, 
coeff$(x_1x_2 \dots x_k)=t-1$. For the case when $k>r$, we generalize Claim 1 to the following: 

{\bf{Claim 2:}} 

For $k>r$, coeff$(x_1x_2 \dots x_k)=(t-1)$ (coeff$(x_2x_3 \dots x_k)+$ coeff$(x_3x_4...x_k)+...+$ coeff$(x_rx_{r+1}...x_k))$. 
(i.e. coeff$(x_1x_2...x_k)=(t-1)$ (coeff$(x_1x_2...x_{k-1})+$ coeff$(x_1x_2...x_{k-2})+...+$ coeff$(x_{1}x_{2}...x_{k-r+1})$.)

The proof of Claim 2 directly generalizes from the proof of Claim 1. Now one mark contributes a factor of $(t-1)$ instead of $(-1)$ 
to the weight of a marked word. For example, for the clusters associated with the word 123456, and $r=3$, the first mark has to be 123, 
the second mark can be 234 or 345. So coeff$(x_1x_2x_3x_4x_5x_6)=(t-1)($coeff$(x_2x_3x_4x_5x_6)+$coeff$(x_3x_4x_5x_6)$). 
In general, like in the proof of Theorem 1, if we are interested in keeping track
of the number of  appearances of the consecutive pattern $12 \dots r$, then there are $(r-1)$ scenarios
of clusters that can give rise to the monomial $x_1x_2 \dots x_k$, depending on where the second mark is. 
By peeling off the first mark, now we loose a factor of $(t-1)$ instead of $(-1)$ in the coefficient of our monomial. 

As the coefficients of the monomials of the same length are the same, Claim 2 immediately implies that 
$weight(C)=\sum_{k=r}^{n} P_k^{(r)}(t) e_k$ where $P_k^{(r)}(t)$ satisfies the recurrence 
$$
P_k^{(r)}(t) \, = \, (t-1) \sum_{i=1}^{r-1} P_{k-i}^{(r)}(t) \quad .
$$
(In fact $P_k^{(r)}(t)$ is just a concise way of writing coeff$(x_1x_2 \dots x_k)$,
 where the consecutive pattern of interest is $12 \dots r$.) From this Theorem 2  follows directly.
\bigskip
{\bf Acknowledgment}: Many thanks are due to Sergi Elizalde for help with the references. Also many thanks to Justin Troyka for pointing out
that ``our" Theorem 1 appeared in Ira Gessel's  PhD thesis.

\bigskip

{\bf References}

[BaNZ]  Andrew Baxter, Brian Nakamura, and Doron Zeilberger, {\it Automatic generation of theorems and proofs on enumerating consecutive-Wilf classes},
in: {\it ``Advances in Combinatorics: Waterloo Workshop in Computer Algebra, W80, May 26-29, 2011"} [Volume in honor of Herbert S. Wilf], 
edited by Ilias Kotsireas and Eugene Zima, 121-138 \quad .

[Bu] Alex Burstein, 
{\it ``Enumeration of words with forbidden patterns"}, Ph.D. thesis, University of Pennsylvania, 1998.

[BuM] Alex Burstein and Toufik Mansour, {\it Words restricted by patterns with at most 2 distinct letters},
The electronic journal of combinatorics 9(2) (2002), \#R3 \quad .

[DB] F. N. Davis and D. E. Barton, {``Combinatorial Chance''}, Hafner, New York, 1962.

[EN] Sergi Elizalde and Marc Noy, {\it Consecutive patterns in permutations},
Advances in Applied Mathematics {\bf 30} (2003), 110-125 .

[Ge]  Ira M. Gessel, {\it ``Generating Functions and Enumeration of Sequences''}, PhD thesis, Massachusetts Institute of Technology, 1977.

[GJ] Ian Goulden and David Jackson, {\it An inversion theorem for cluster decomposition of sequences with distinguished 
subsequences},  J. London Math. Soc.(2) {\bf 20} (1979), 567-576.

[MR] A. Mendes and J.B. Remmel, {\it Permutations and words counted by consecutive patterns}, Adv. Appl. Math, {\bf 37}(2006), 443-480.

[Na] Brian Nakamura, {\it Computational approaches to consecutive pattern avoidance in permutations},
Pure Math. Appl. (PU.M.A.) {\bf 22}(2011), 253-268. \hfill\break
Also available from {\tt https://arxiv.org/abs/1102.2480} \quad .

[NoZ] John Noonan and Doron Zeilberger, {\it The Goulden-Jackson cluster method: extensions, applications, and implementations},
J. Difference Eq. Appl. {\bf 5}(1999), 355-377.  \hfill\break
{\tt http://sites.math.rutgers.edu/\~{}zeilberg/mamarim/mamarimhtml/gj.html} \quad .

[P] Lara Pudwell, {\it Enumeration schemes for pattern-avoiding words and permutations}, Ph.D. thesis, Rutgers University, May 2008, \hfill\break
{\tt http://faculty.valpo.edu/lpudwell/papers/pudwell\_thesis.pdf} \quad .

[W] Xiangdong Wen (2005) {\it The symbolic Goulden-Jackson cluster method}, Journal of Difference Equations and Applications, 11:2, 173-179, DOI: 10.1080/10236190512331329432 \hfill\break
{\tt 
https://doi.org/10.1080/10236190512331329432} \quad .

[Wi] The Wikipedia Foundation, {\it Permutation Pattern} [Initiated and largely written by Vince Vatter] \hfill\break
{\tt https://en.wikipedia.org/wiki/Permutation\_pattern} \quad .

\bigskip
\hrule
\bigskip
Mingjia Yang, Department of Mathematics, Rutgers University (New Brunswick), Hill Center-Busch Campus, 110 Frelinghuysen
Rd., Piscataway, NJ 08854-8019, USA. \hfill\break
Email: {\tt my237 at math dot rutgers dot edu}   \quad .
\bigskip
Doron Zeilberger, Department of Mathematics, Rutgers University (New Brunswick), Hill Center-Busch Campus, 110 Frelinghuysen
Rd., Piscataway, NJ 08854-8019, USA. \hfill\break
Email: {\tt DoronZeil at gmail  dot com}   \quad .
\bigskip
\hrule
\bigskip

First written: May 15, 2018; This version: May 22, 2018.

\end